\documentclass[review,a4paper]{elsarticle}
\usepackage{amsmath}
\usepackage{color}
\usepackage{graphicx}
\usepackage{bm}
\usepackage{amssymb}
\usepackage[utf8]{inputenc}

\usepackage{lineno,hyperref}
\modulolinenumbers[5]

\usepackage{setspace}
\usepackage[section]{placeins}

\journal{Chaos, Solitons \& Fractals}

\biboptions{numbers,sort&compress}

\begin{document}

\begin{frontmatter}

\title{A filtered Hénon map}

\author{Vinícius S. Borges}
\corref{mycorrespondingauthor}
\cortext[mycorrespondingauthor]{Corresponding author.}
\ead{viniciusb@usp.br}
\author{Marcio Eisencraft}
\ead{marcioft@usp.br}
\address{Telecommunication and Control Engineering Department, Escola Politécnica, \\University of São Paulo, Brazil\\Published in Chaos, Solitons $\&$ Fractals \url{https://doi.org/10.1016/j.chaos.2022.112865}}

\begin{abstract}
In this paper, we use Lyapunov exponents to analyze how the dynamical properties of the Hénon map change as a function of the coefficients of a linear filter inserted in its feedback loop. We show that the generated orbits can be chaotic or not, depending on the filter coefficients. The dynamics of the system  presents complex behavior, including cascades of bifurcations, coexistence of attractors, crises, and ``shrimps".  The obtained results are relevant in the context of bandlimited chaos-based communication systems, that have recently been proposed in the literature.
\end{abstract}


\end{frontmatter}


\section{Introduction}
\label{ch::introducao}

A chaotic signal has three main characteristics: it is bounded, presents aperiodicity and sensitive dependence on the initial conditions (SDIC) \cite{Alligood2000}. These properties have stimulated proposals of chaotic signals implementation proposals in Telecommunications and Signal Processing since the seminal Pecora and Carroll work \cite{Pecora1990}. They have shown that two identical systems generating chaotic signals could be synchronized despite SDIC. Since then, many possible applications such as chaos-based communication systems (CBCS) \cite{Souza2019,Baptista2021}, watermarking \cite{Loan2018}, compressed-sensing \cite{Rontani2016}, image encryption \cite{Li2015,Zhou2016}, ultra-wideband communications \cite{Dmitriev2006}, memristor models \cite{peng2020discrete} and others have appeared.

Recently, the performance of CBCS has been accessed in real transmission scenarios involving channel distortion, noise, bandwidth constraints and delay \cite{Liu2020,Baptista2021}. Since transmission channels are always bandlimited \cite{Lathi2009}, it is necessary to know and control the bandwidth of the transmitted chaotic signals. In this regard, the authors of \cite{Eisencraft2009b} employed a discrete-time non recursive linear filter built into the chaos generator as a way to control the bandwidth of chaotic signals. After that, it was demonstrated that the filter insertion does not affect chaotic synchronization \cite{Fontes2016a}, which is essential for chaos communication. However, the question remained as to under what conditions the generated signals were still chaotic. The chaotic nature of the transmitted signals issue was only barely touched in \cite{Eisencraft2011a}.

Bearing this in mind,
in the present paper, we analyze the dynamics of the discrete-time dynamical system obtained when we add a two coefficients non recursive linear filter in the feedback loop of the Hénon map \cite{henon1976two}.This system is a simplified version of the one considered in \cite{Eisencraft2009b,Fontes2016a}. Using only two coefficients permits a more insightful analysis using dynamical systems tools. Besides, non recursive filters are always BIBO (bounded input, bounded output stable) stable \cite{Oppenheim2009}, so any divergence presented by the orbits is caused by the dynamics and not by the filter itself.

This paper is organized as follows: the system under consideration is described in Section 2. Its dynamical analysis and the main results are presented in Section 3 and our conclusions are drafted in Section 4.

\section {The filtered Hénon map}
\label{ch::filtered_henon}

The Hénon map is given by \cite{henon1976two}
\begin{equation}
\begin{cases}
x_{1}(n+1) = \alpha - \left(x_{1}(n)\right)^2 + \beta x_{2}(n) \\
x_{2}(n+1) = x_{1}(n)
\end{cases}
\label{eq:henon}
\end{equation}
where ${\alpha, \beta}$ are real parameters, $n=0,1,\ldots$, and the state variables are \linebreak$\bm x=\left[x_1,x_2\right]$. It has been used as a paradigm to generate 2-dimensional discrete-time chaotic signals \cite{lellep2020using,zhao2020efficient}.

In \cite{Eisencraft2009b} it was proposed to filter $x_1(n)$ using a non recursive or finite impulse response (FIR) filter \cite{Oppenheim2009} so that
\begin{equation}
x_{3}(n)=\sum_{j=0}^{N_{S}-1} c_{j} x_{1}(n-j),
\label{eq:FIR}
\end{equation}
where $c_{j}, 0 \leq j \leq N_{S}-1$ are the filter coefficients. Given $N_{S}$ and a frequency response specification, there is a number of different FIR filter design techniques to obtain the coefficients $c_{j}$ \cite{Oppenheim2009}.

To allow synchronization in the receiver, $x_3(n)$ is fed back in the dynamical system, so that the resulting system is
\begin{equation}
\begin{cases}
x_{1}(n+1)=\alpha-\left(x_{3}(n)\right)^{2}+\beta x_{2}(n) \\
x_{2}(n+1)=x_{1}(n) \\
x_{3}(n+1)=\sum_{j=0}^{N_{S}-1}{ c_{j} x_{1}(n-j{\color{blue}+1})}
\end{cases}.
\label{eq:henon_FIR}
\end{equation}

In \cite{Fontes2016a} this system was used as a way to generate a chaotic low-pass signal so that it could be transmitted through a communication channel without distortion. In addition, the authors were able to prove that synchronization is achieved independently of the filter coefficients. When it comes to the chaotic nature of the orbits of \eqref{eq:henon_FIR}, preliminary numerical results show that depending on $c_j$, the orbits can still be chaotic, can converge to periodic attractors or diverge towards infinity \cite{Eisencraft2011a}.

Here we consider the simplest non-trivial case of $N_S=2$ coefficients. In this case,
\begin{equation}
x_{3}(n+1)=c_{0} x_{1}(n+1)+c_{1} x_{1}(n)
\end{equation}
and so \eqref{eq:henon_FIR} becomes
\begin{equation}
\begin{cases}
x_{1}(n+1)=\alpha-\left(x_{3}(n)\right)^{2}+\beta x_{2}(n) \\
x_{2}(n+1)=x_{1}(n) \\
x_{3}(n+1)=c_{0} x_{1}(n+1)+c_{1} x_{1}(n)
\end{cases}
\label{eq:henon_FIR_2d}
\end{equation}
or
\begin{equation}
\begin{cases}
x_{1}(n+1)=\alpha-(c_{0} x_{1}(n)+c_{1} x_{2}(n))^{2}+\beta x_{2}(n) \\
x_{2}(n+1)=x_{1}(n).
\label{eq:henon_FIR_2}
\end{cases}
\end{equation}
Thus, the resulting system is a 2-dimensional non-linear dynamical system that has the original Hénon map \eqref{eq:henon} as a particular case for $c_0=1$ and $c_1=0$. In what follows we unravel the dynamical properties of this system as a function of the filter parameters $c_0$ and $c_1$, considering $\alpha = 1.4$ and $\beta = 0.3$ fixed as in \cite{henon1976two}.

\section{Analysis of the map dynamics}

We begin studying the fixed points of \eqref{eq:henon_FIR_2} and their stability. Then, we perform a numerical analysis of the periodicity of its orbits and of the largest Lyapunov exponent of the attractors.

\subsection{Fixed Points}

From \eqref{eq:henon_FIR_2}, we can calculate the fixed points $\bm{P}_i=\left(p_i,p_i\right)$, $i=1,2$ of the filtered Hénon map solving
\begin{equation}
p_i = \alpha-\left(c_{0} p_i+c_{1} p_i\right)^{2}+\beta p_i\Rightarrow (c_{0}+c_{1})^{2}p_i^{2}+(1-\beta)p_i-\alpha=0.
\end{equation}

We obtain
\begin{align}
p_{1} &= \frac{-(1-\beta)- \sqrt{(1-\beta)^{2} + 4\alpha (c_{0} +c_{1})^{2}}}{2(c_{0} +c_{1})^{2}}<0
\label{eq::pontofixop1}
\end{align}
and
\begin{align}
p_{2} &= \frac{-(1-\beta)+ \sqrt{(1-\beta)^{2} + 4\alpha (c_{0}+c_{1})^{2}}}{2(c_{0} +c_{1} )^{2}}>0.
\label{eq::pontofixop2}
\end{align}
This way, except for the case $c_0=c_1=0$, we have two fixed points $\bm{P}_1=\left(p_1,p_1\right)$ and $\bm{P}_2=\left(p_2,p_2\right)$ that depends only on the squared sum of the filter coefficients $\left(c_0+c_1\right)^2$. For $c_0=c_1=0$ the system has only one fixed point at $\left(\frac{\alpha}{1-\beta},\frac{\alpha}{1-\beta}\right)$.

Stabilities of $\bm{P}_1$ and  $\bm{P}_2$ are determined by the largest absolute value of the eigenvalues of the Jacobian matrix of \eqref{eq:henon_FIR_2} \begin{equation}
\bm{J}=
\begin{bmatrix}
-2c_0(c_0x_1+c_1x_2)  & - 2c_{1}(c_0x_1+c_1x_2)+\beta \\
1 & 0
\end{bmatrix},
\end{equation}
calculated on $\bm{P}_1$ and $\bm{P}_2$, respectively. They are given by
\begin{align}
	\lambda= \max&\left\{\left|-pc_0\left(c_0+c_1\right)+\sqrt{\left(pc_0\left(c_0+c_1\right)\right)^2-2pc_1\left(c_0+c_1\right)+\beta}\right|,\right.\nonumber\\&\left.\left|-pc_0\left(c_0+c_1\right)-\sqrt{\left(pc_0\left(c_0+c_1\right)\right)^2-2pc_1\left(c_0+c_1\right)+\beta}\right|\right\}
	\label{eq:eigvalues_1},
\end{align}
where $p$ should be substitute for $p_1$ of \eqref{eq::pontofixop1} for $\bm{P_1}$ or $p_2$ of \eqref{eq::pontofixop2} for $\bm{P_2}$. If $\lambda>1$ we have an unstable fixed point and when $\lambda<1$ we have a stable fixed point \cite{Alligood2000}.

\begin{figure}[htb]
\centering
\includegraphics[width=\textwidth]{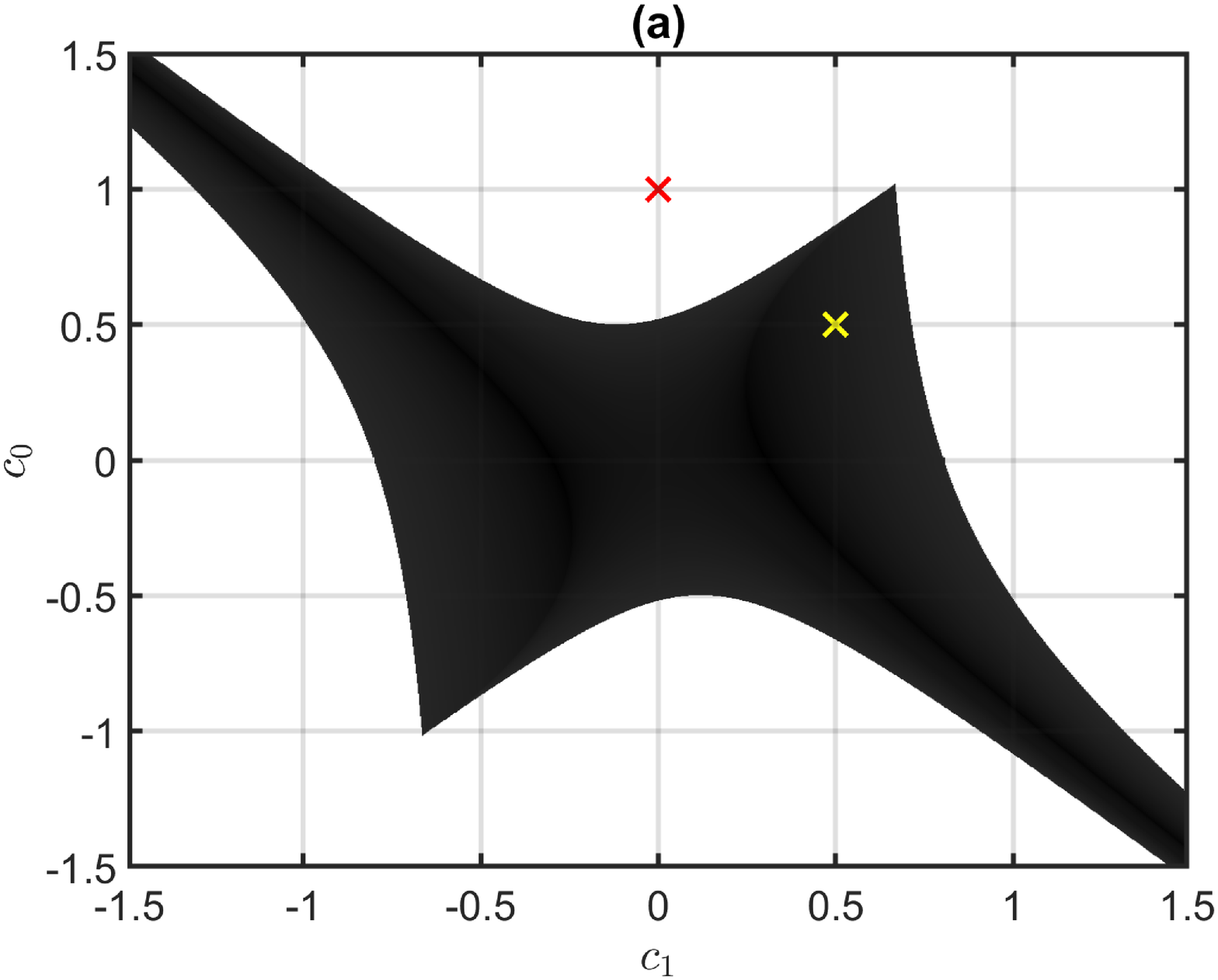}
\includegraphics[width=0.49\textwidth]{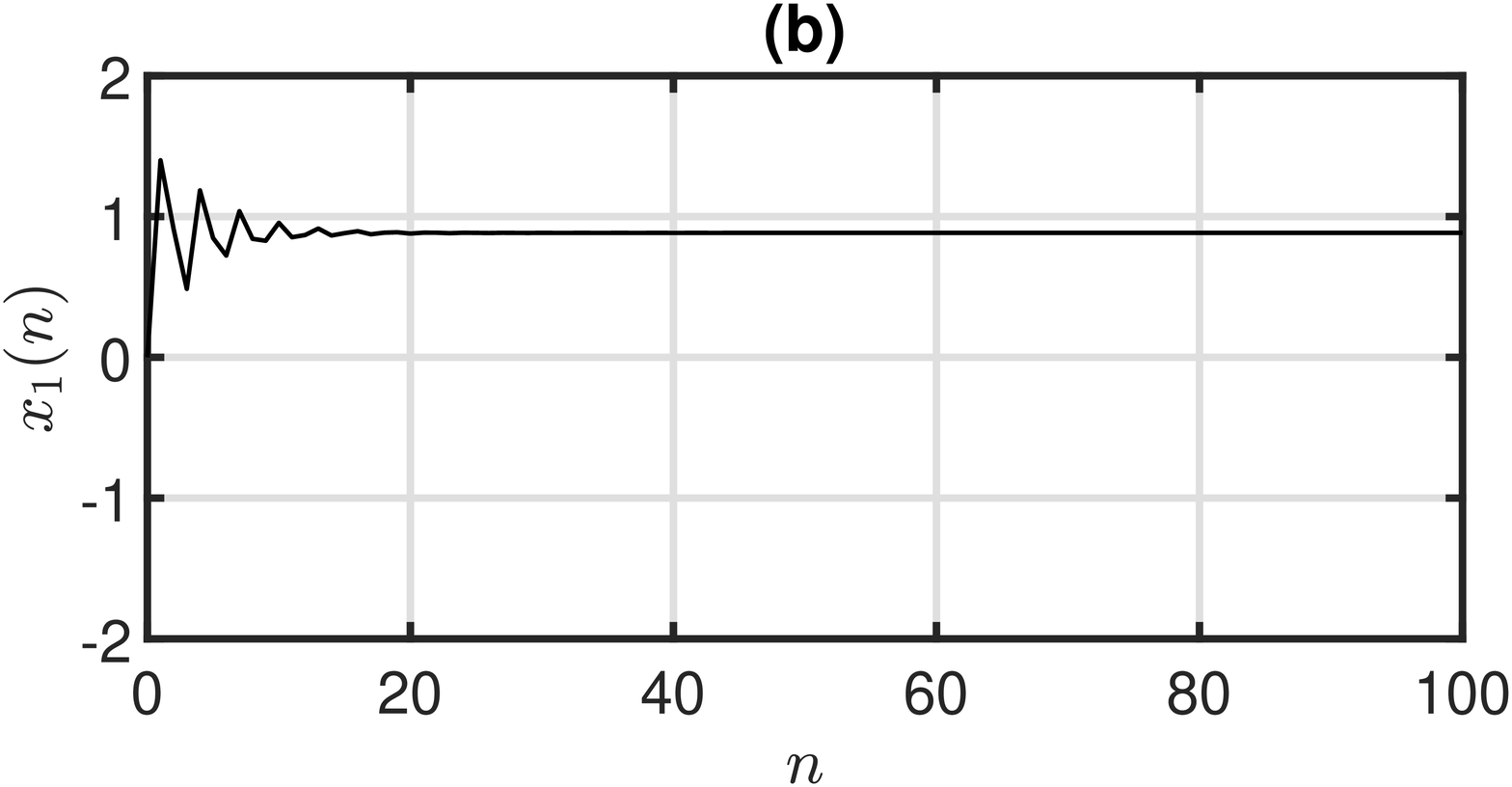}
\includegraphics[width=0.49\textwidth]{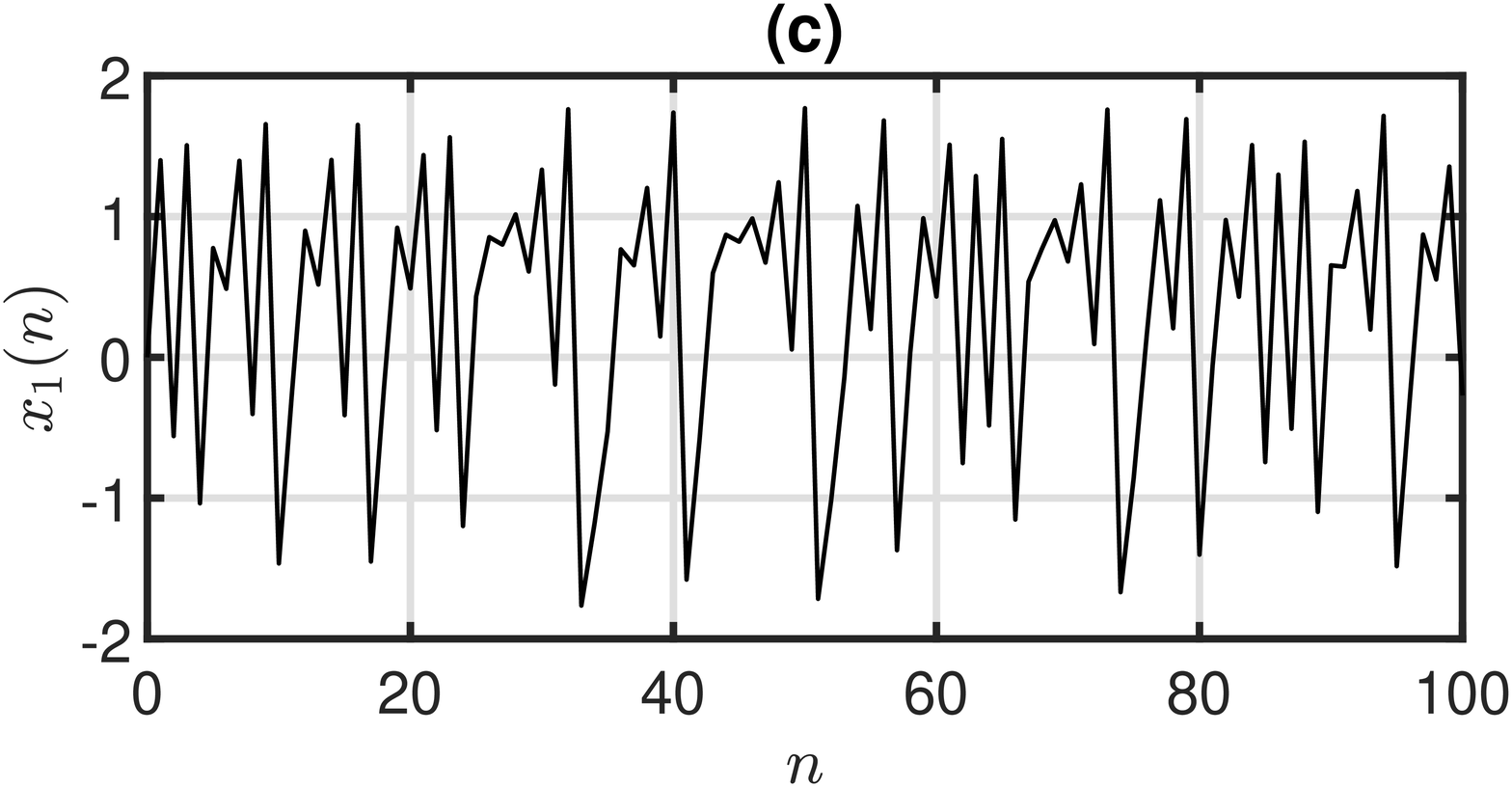}
 \caption{(a) Stability region of $\bm{P}_2$ is shown in black in the $c_0\times c_1$ plane; examples of $x_1(n)$  for (b) $c_{0}=c_{1}=0.5$ (yellow cross in (a)) leading to a stable $\bm{P}_2$ and (c) $c_{0}=1$, $c_{1}=0$ (red cross in (a)) leading to an unstable $\bm{P}_2$ of the original Hénon map, i.e, without filter.}
\label{fig:pontos_fixos}
\end{figure}

For $\bm{P_1}$ we numerically found that $\lambda>1$ for all values of $c_0$ and $c_1$  so that it is always unstable. When it comes to
$\bm{P_2}$, there is a region of stability in the $c_{0}\times c_{1}$ plane, where $\lambda<1$. This region is shown in black in Figure \ref{fig:pontos_fixos}(a). It presents an odd symmetry with respect to the axis $c_1=0$ and it is unbounded in the direction of the bisector of the 2st and 4rd quadrants $c_0+c_1=0$. Typical examples of $x_1(n)$ for parameters leading to stable and unstable $\bm{P_2}$ are shown in Figure \ref{fig:pontos_fixos}(b) and (c), respectively. In case (c), we considered $c_0=1$ and $c_1=0$, which leads to the original chaotic Hénon map \cite{henon1976two}.

\subsection{Largest Lyapunov Exponent}

The $k$-th Lyapunov exponent for an orbit with initial condition $\bm{x}(0)$ is defined as
\begin{equation}
	h_{k}\left(\bm{x}_{0}\right)=\lim _{n \rightarrow \infty}\frac{1}{n} \ln \left(\left\|\bm{J}^{n}\left(\bm{x}(0)\right) \bm{u}_{\bm{k}}\right\|\right),
	\label{eq::lyap}
\end{equation}
where $\bm{J}^{n}$ is the Jacobian matrix of the $n$-time iterated map evaluated at point $\bm{x}(0)$ and $\bm{u}_{\bm{k}}$ is the eigenvector corresponding to the $k$-th largest eigenvalue of this Jacobian matrix. For sake of simplicity of notation we define $h_0\triangleq h$.

In order to verify the presence of chaotic orbits generated by \eqref{eq:henon_FIR_2d}, we numerically analyze the largest Lyapunov  exponent  $h$ of its orbits. Chaotic behavior can be identified in a bounded aperiodic orbit when $h>0$ \cite{Alligood2000}. Numerical estimators of \eqref{eq::lyap} can be obtained by a variety of techniques. Here,  we consider the tangent map method \cite{Alligood2000}.

For the numerical evaluations of $h$, we considered random initial conditions uniformly distributed in the unit square and 3000 iterations, excluding the first  500 iterations. For each estimate, we took 25 different initial conditions, and present the mean value. Our numerical experiments have shown that these quantities of initial conditions and iterations are sufficient to determine $h$ of the attractor with an accuracy higher than $10^{-2}$.

\begin{figure}[htb]
 	\centering
\includegraphics[width=\textwidth]{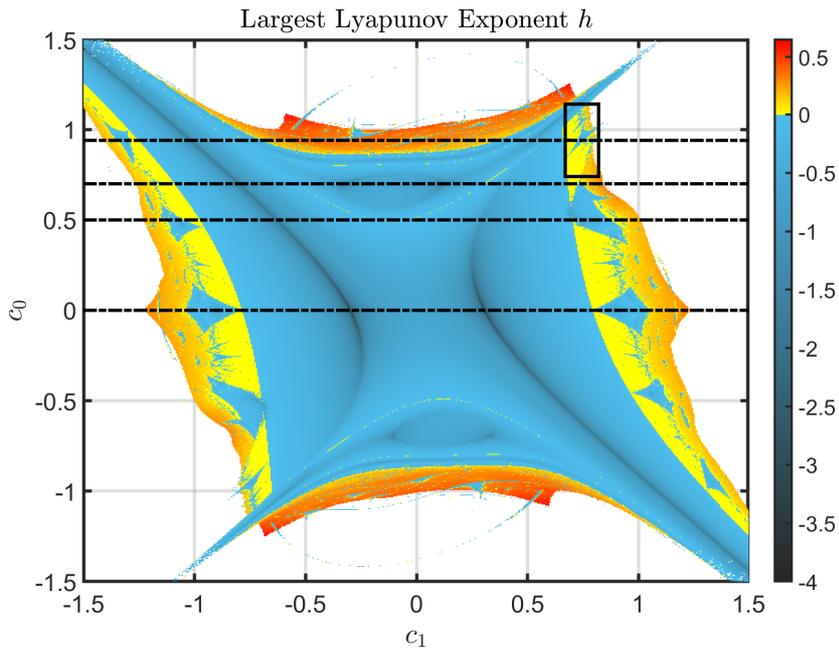}
\caption{Largest Lyapunov exponent $h$ for the map \eqref{eq:henon_FIR_2} as a function of $c_0$ and $c_1$. White regions represent parameters that generate divergent orbits. The region marked with a black rectangle is further analyzed in Figure \ref{fig:Henonfilterzoon1}}.
\label{fig:Henonbifu1}
\end{figure}

Figure \ref{fig:Henonbifu1} shows a general picture of how $h$ varies in the plane $c_{0} \times c_{1}$ inside the square $[-1.5,1.5]\times[-1.5,1.5]$.

\begin{figure}[htb]
\centering
\includegraphics[width=1\textwidth]{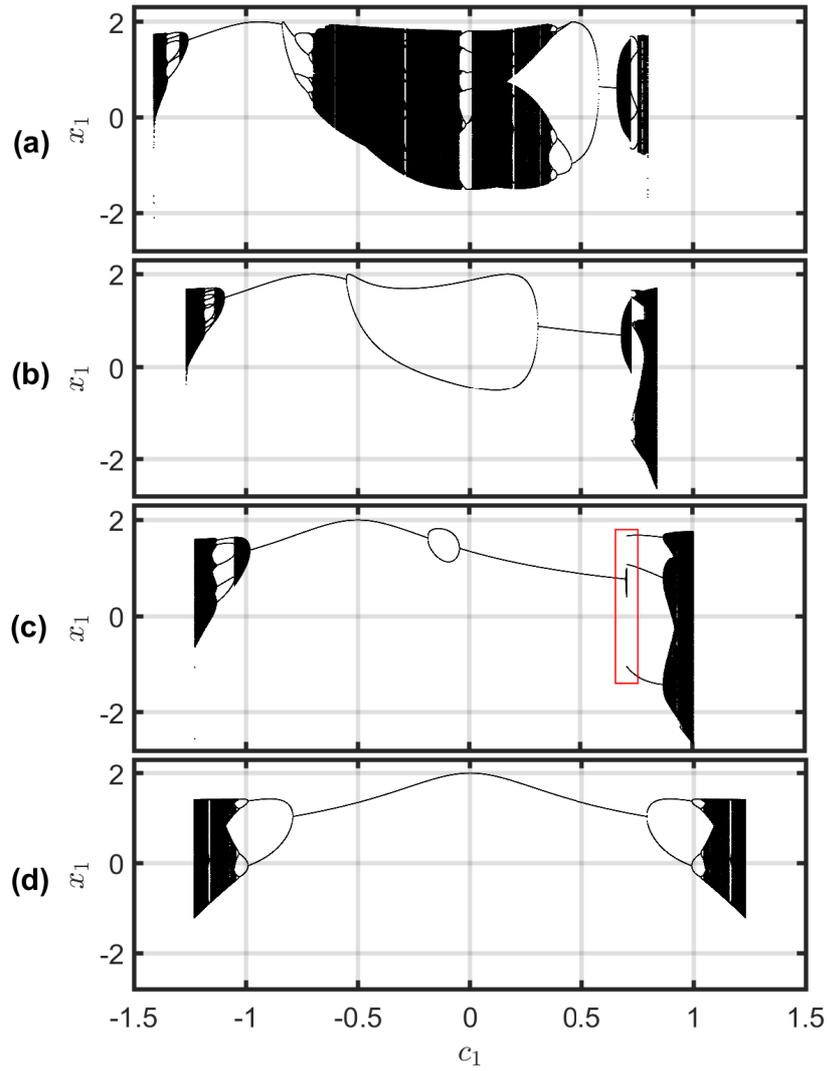}
\vspace{-1.5cm}
\caption{Bifurcation diagram as a function of $c_1$ for (a) $c_{0}=0.94$,  (b) $c_{0}=0.70$, (c) $c_{0}=0.50$ and (d) $c_{0}=0$. These values of $c_0$ are marked by dashed lines in Figure \ref{fig:Henonbifu1}.}
 	\label{fig:Henonbifu2}
 \end{figure}

We can clearly see a central region where $h<0$,  shown in different shades of blue. As expected, this region encompasses the black one in Figure \ref{fig:pontos_fixos}(a) where there is a stable fixed point. This central region contains orbits that converges towards $\bm{P_2}$ and period-2 orbits, as can be seen in the bifurcation diagrams of Figure \ref{fig:Henonbifu2}. They were plotted as a function of $c_1$ for constant values of $c_0$, indicated by the dashed lines in Figure \ref{fig:Henonbifu1}. The strategy of following the attractor was applied, using as initial condition for each value of $c_1$ a point of the attrator obtained for the previous value of $c_1$.

Surrounding the blue central region in Figure \ref{fig:Henonbifu1}, there are ``chaotic regions'' where $h>0$ predominates (in shades of yellow and red) interspersed by islands where $h<0$ but with higher periodicity, as the biffurcation diagrams of Figure \ref{fig:Henonbifu2} shows.

Increasing even more $c_0$ and $c_1$ results in orbits that diverge towards infinity, except for the region in the vicinity of the bisector of the 2st and 4rd quadrants $c_0+c_1=0$. These divergence regions are shown in white in Figure \ref{fig:Henonbifu1} and by the regions without attractor points in the bifurcation diagrams of Figure \ref{fig:Henonbifu2}.

In Figure \ref{fig:Henonbifu2}(a), one can see that for $c_{0}=0.94$,  the transitions from the blue region to the ``chaotic regions'' sometimes present period-doubling cascades and sometimes are abrupt, via crisis. Inside the chaotic regions there are windows of periodicity associated with the blue islands in their inside. The transitions from chaos to divergence are abrupt in both ends of the diagram.

Figure \ref{fig:Henonbifu2}(b) shows a large region with a period-2 stable orbit inside the blue region of Figure \ref{fig:Henonbifu1} for $c_0=0.7$. The transition to chaos is abrupt in the end of the fixed point $\bm{P_2}$ stability region.

\begin{figure}[htb]
	\centering
	\includegraphics[width=1\columnwidth]{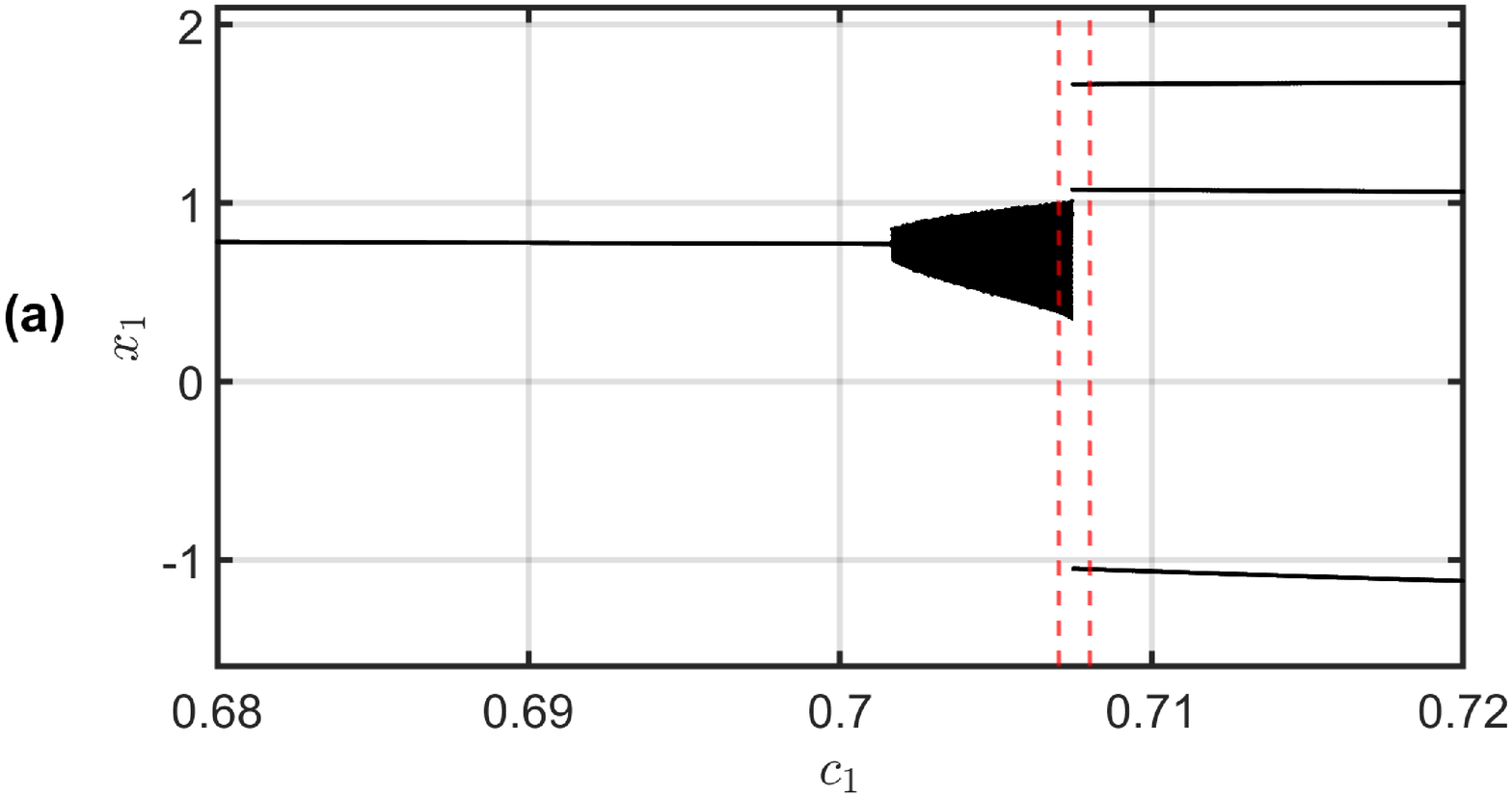}
	\includegraphics[width=\columnwidth]{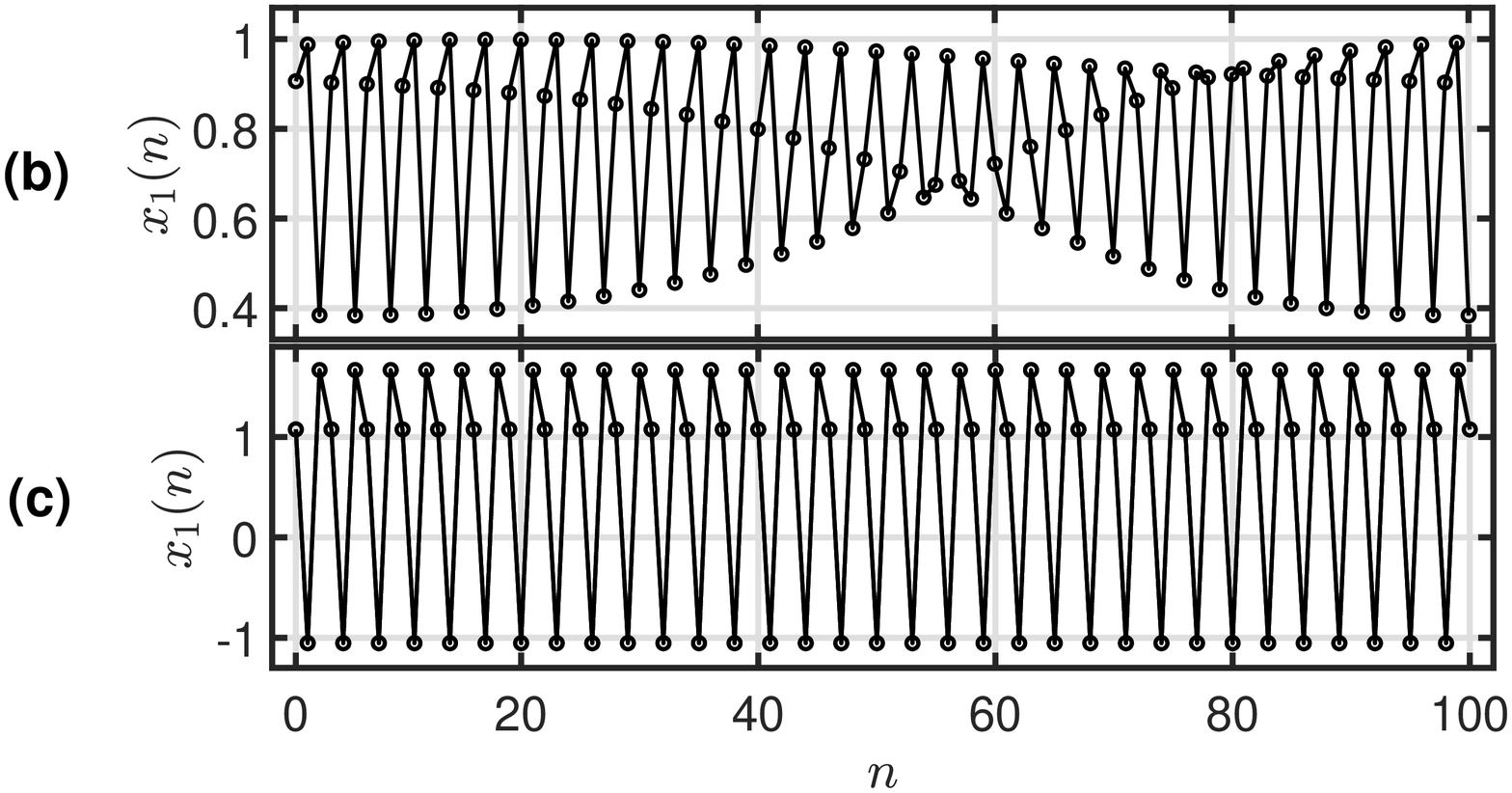}
    \caption{(a) Zoom of the region indicated by the red rectangle in Figure \ref{fig:Henonbifu2}(c); example of orbits for (b) $c_1=0.707$ and initial condition $\bm x(0)=[0.8,0.8]$; (c) $c_1=0.708$. These values of $c_1$ are marked by red dashed lines in (a).}
	\label{fig:Henonorbit1}
\end{figure}

A similar patter is presented in Figure \ref{fig:Henonbifu2}(c) for $c_{0}=0.5$, with chaos or quasi-periodicity appearing and disappearing abruptly with changing $c_1$. As an example of this behavior, Figure \ref{fig:Henonorbit1}(a) presents a zoom for $c_{1}\in[0.68;0.72]$ (region marked with a red rectangle in Figure \ref{fig:Henonbifu2}(c)). Figures \ref{fig:Henonorbit1}(b) and (c) show examples of an aperiodic orbit for $c_1=0.707$ and a period-3 orbit for $c_1=0.708$ respectively. These values of $c_1$ are signaled by red dashed lines in Figures \ref{fig:Henonorbit1}(a).

It is relevant to note that in the region of Figure \ref{fig:Henonorbit1} where there are aperiodic orbits, we observed the presence of two coexistent attractors with large basins of attraction. For instance, in $c_1=0.707$ the aperiodic orbit with $h=0.00052\pm0.00044$ is obtained for $\left\{x_1(0),x_2(0)\right\}\subset[0.7,0.9]$ and a period-3 orbit with $h=-0.4764\pm0.0020$ is obtained for $\left\{x_1(0),x_2(0)\right\}\subset[0,0.7]$.

In figure \ref{fig:Henonbifu2}(d), the diagram was obtained for $c_0=0$. In this case, the part of the cutting in the $c_1\times c_2$ plane that crosses the blue region in Figure \ref{fig:Henonbifu1} is completely contained in the stability region of $\bm{P_2}$ (black region in Figure \ref{fig:pontos_fixos}). This way, there is no period-2 orbits and the fixed point $\bm{P_2}$ is present along all this part. When the values of $c_1$ passes the blue central region of Figure \ref{fig:Henonbifu1}, we have a period-doubling cascade transition to chaos.

\begin{figure}[htb]
	\centering
\includegraphics[width=\linewidth]{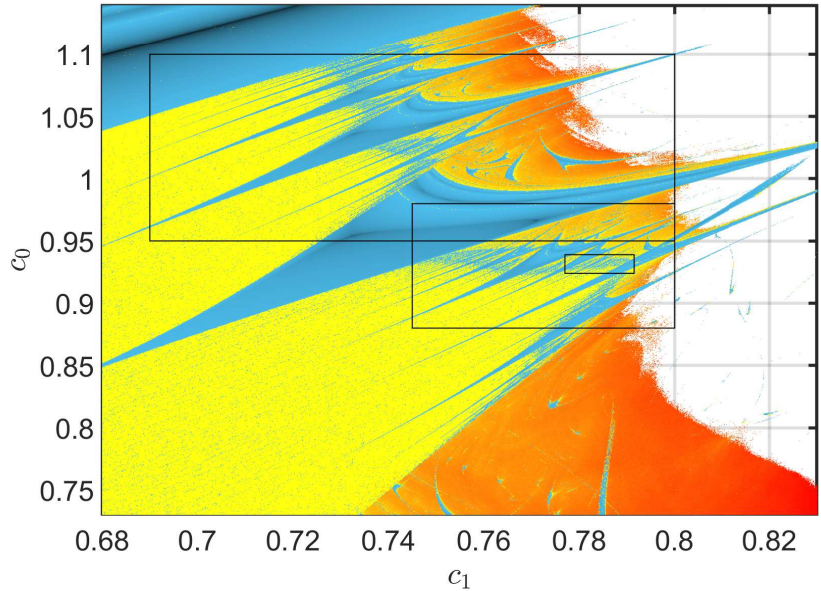}
	\caption{Zoom of the region $0.74 \leq c_{0} \leq 1.14$, $0.68 \leq c_{1} \leq 0.83$ of Figure \ref{fig:Henonbifu1}. The regions within the indicated rectangles are analyzed in more detail in Figures \ref{fig:Henonfilterzoon2} (largest and medium rectangles) and \ref{fig:Henonfilterzoon3} (smallest rectangle).
\label{fig:Henonfilterzoon1}}
\end{figure}

\begin{figure}[htb]
\centering
\includegraphics[width=\linewidth]{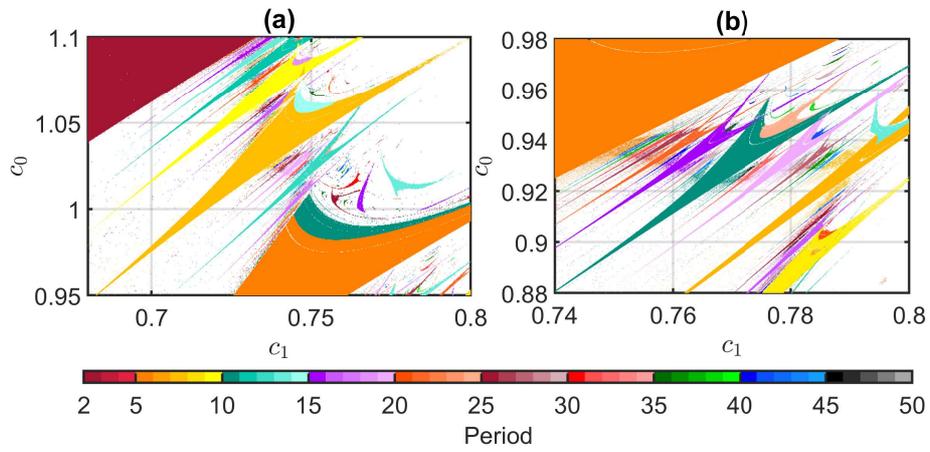}
\caption{Zooms of the (a) larger and (b) medium rectangles in Figure \ref{fig:Henonfilterzoon1} highlighting the presence of different periodic shrimps. Different colors represent different periods as shown in the color bar. White regions represents chaos or divergence. \label{fig:Henonfilterzoon2}}
\end{figure}

In Figure \ref{fig:Henonfilterzoon1} we explore with more detail the region $0.74 \leq c_{0} \leq 1.14$, $0.68 \leq c_{1} \leq 0.83$ that is marked with a black rectangle in Figure \ref{fig:Henonbifu1}. One can clearly notice fractal-like patterns of periodic islands, popularly known as \emph{shrimps} \cite{de2012self,gallas1993structure,dos2016unstable}, surrounded by chaotic regions. Each shrimps presents orbits with different periodicity, as can be seen in the enlargements shown in Figure \ref{fig:Henonfilterzoon2}.

\begin{figure}[htb]
\centering
 	\includegraphics[width=.8\linewidth]{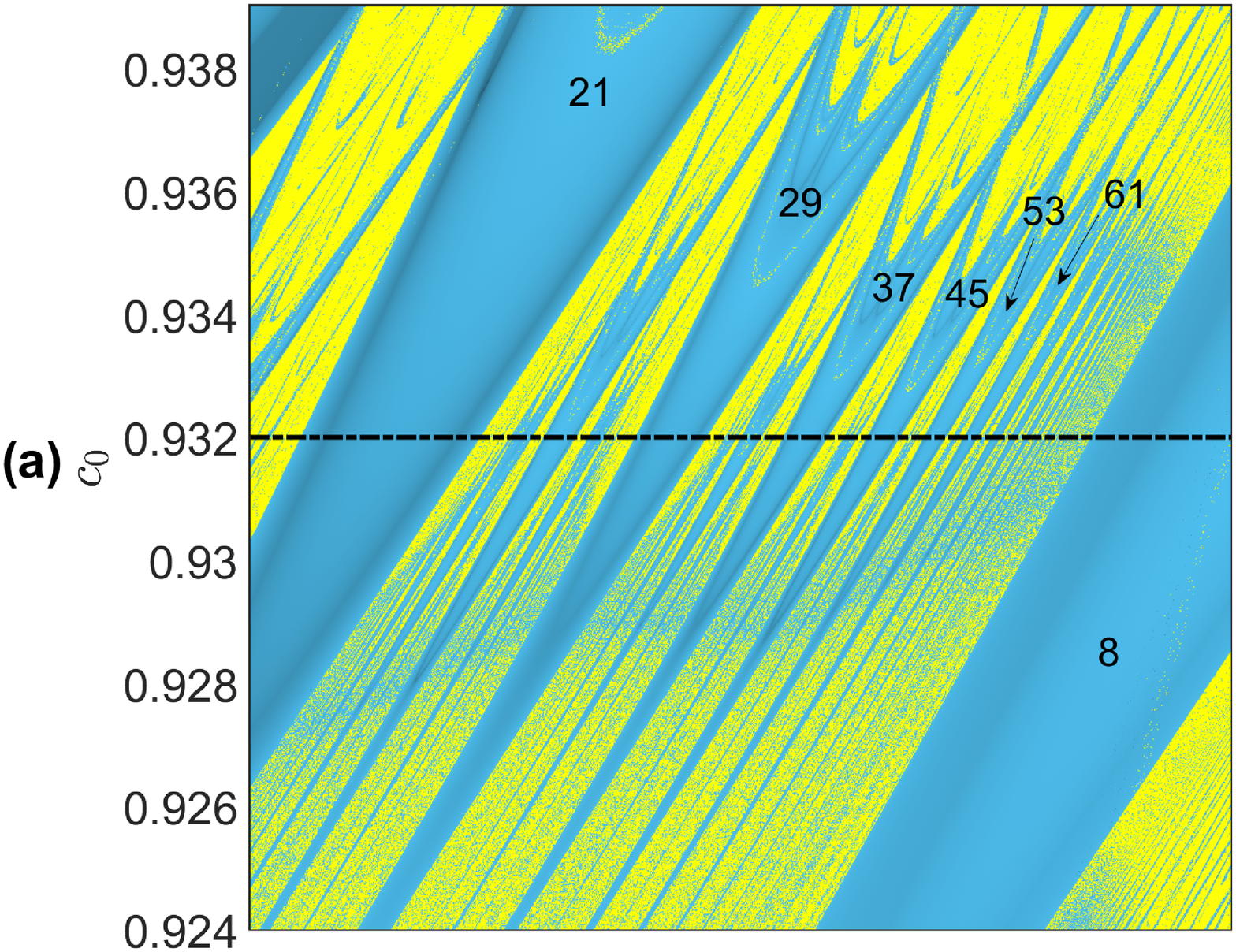}
 	\includegraphics[width=0.71\linewidth]{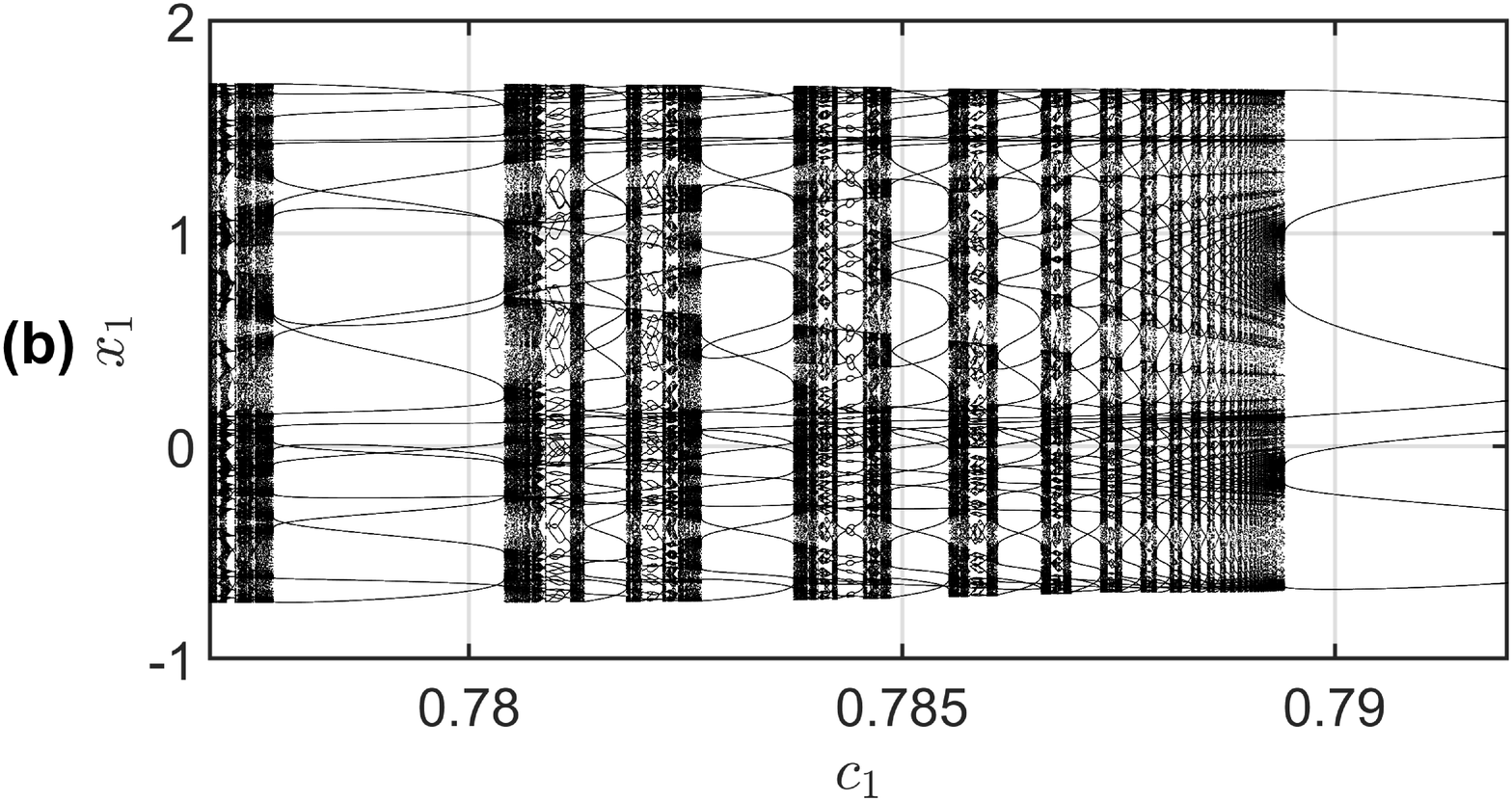}
  	\caption{(a) Zoom of the smallest rectangle in Figure \ref{fig:Henonfilterzoon1}. Periods of the orbits inside some shrimps are indicated; (b) Bifurcation diagram for $c_0=0.932$ (indicated by a dashed line in (a)).\label{fig:Henonfilterzoon3}}
\end{figure}

The fractal nature of the distribution of shrimps and their different periods are highlighted in the zoom of the smallest rectangle in Figure \ref{fig:Henonfilterzoon1} presented in Figure \ref{fig:Henonfilterzoon3}(a). In Figure \ref{fig:Henonfilterzoon3}(b) a bifurcation diagram for $c_0=0.932$ (indicated by the dashed line in Figure \ref{fig:Henonfilterzoon3} (a)) is show. The presence of multiple periodic windows of high period and a cascading of windows is clearly visible. The chaotic structure appears immediately when crossing the periodic shrimps boundaries without a period-doubling or other classical bifurcation scenario.

The presence of shrimps, especially their thin antennae, in the parameter region in which chaotic orbits are generated is a  challenge from the point of view of practical applications of chaotic signals.  Small perturbations in the values of the chosen coefficients can cause the chaotic behavior to be replaced by a stable periodic one, losing the DSCI that is critical in applications involving secure CBCS \cite{Shao,Vaseghi,wang2004new}. This issue seems to have gone unnoticed in many previous papers employing map-generated chaos in communications \cite{Eisencraft2012aa,Eisencraft2009b,Fontes2016a}, and analyses like the ones presented here need to be taken into account when choosing the map and filters to be used to generate bandlimited chaotic signals.

\section{Conclusions}

In this paper we present an analysis of the Hénon map including a linear non recursive filter with two coefficients in the feedback loop. The use of such filters is a way of generating bandlimited chaotic signals for use in CBCS.

Our numerical simulations, by means of the largest Lyapunov exponent, show that the filter coefficients change the chaotic properties of the orbits in a complex way, including the appearance of shrimps. The presence of these structures shows that small changes in the filter coefficients, specially in the parameter space in the vicinity of their thin antennae can completely change the dynamic properties of the generated orbits. They can become chaotic or periodic due to quantization errors expected in any implementation, for example. This way, bandlimited CBCS must be carefully projected to guarantee that the generated signals remain chaotic.

A natural next step is an analysis of what happen if filters with more coefficients or with a recurrent structure are employed , as is usual in communication systems. One possibility is to study how the order and cut-off frequency of low-pass filters generated by different design techniques influence the Lyapunov exponent.

\section*{Acknowledgments}
The authors thank Prof. Antonio M. Batista for interesting discussions throughout the writing of this paper.

This study was financed in part by CNPq-Brazil (grants 140081/2022-4 and 311039/2019-7) and the CAPES-Brazil (Finance Code 001).

\end{document}